\documentclass[twoside,bezier,reqno]{amsart}

\usepackage{epsfig}
\usepackage{amsmath,amsthm,amsfonts,amssymb,amscd,euscript}
\usepackage{graphicx}
\usepackage{savesym}
\usepackage{graphics,amsopn,amstext}
\usepackage{amsrefs}
\usepackage{mathtools}
\usepackage{hyperref}
\hypersetup{
	colorlinks,
	citecolor=blue,
	filecolor=black,
	linkcolor=black,
	urlcolor=black
}
\usepackage{mathptmx}
\usepackage[11pt]{moresize}
\usepackage{fancyhdr}
\usepackage{psfrag}
\usepackage{latexsym}
\usepackage{url}
\usepackage{color}
\usepackage{txfonts}
\usepackage{float}
\usepackage{cases}
\usepackage{algorithmic}
\usepackage{algorithm}
\usepackage{verbatim}
\usepackage{listings}
\usepackage{tikz}
\usetikzlibrary{arrows}
\usetikzlibrary{decorations.pathreplacing}
\usepackage{framed}
\usepackage{cite}
\usepackage{wasysym}
\usepackage{multirow}
\usepackage{xr}
\usepackage{subfig,tikz}
\usepackage{algorithm,algorithmic}
\usepackage[ampersand]{easylist}
\input{epsf}
\newcommand{\filebegin}{\begin{document}}
\newcommand{\fileend}{\end{document}}
\makeatother
\def\thefootnote{}
\newcommand{\lo}{\longrightarrow}
\newcommand{\NMM}{\hspace*{2mm}}

\renewcommand{\baselinestretch}{1.1}

\input{epsf}
\renewcommand{\baselinestretch}{1.1}
\def\n{\noindent}%

\numberwithin{equation}{section}
\def\mapdown#1{\Big\downarrow\rlap
{$\vcenter{\hbox{$\scriptstyle#1$}}$}}
\newtheorem{theorem}{Theorem}[section]
\newtheorem{lemma}[theorem]{Lemma}
\newtheorem{proposition}[theorem]{Proposition}
\newtheorem{corollary}[theorem]{Corollary}

\theoremstyle{definition}
\newtheorem{definition}[theorem]{Definition}
\newtheorem{example}[theorem]{\sc Example}
\newtheorem{xca}[theorem]{Exercise}
\theoremstyle{remark}
\newtheorem{remark}[theorem]{Remark}

\begin{document}

\vspace*{2cm}
\begin{center}
	{\bf\large The Hyper-Zagreb index of trees and unicyclic graphs}
	\\[0.5cm]
	{Hassan Rezapour$^a$, Ramin Nasiri$^{b*}$ \footnote{$^*$Corresponding Author},  Seyedahmad Mousavi$^c$\\[2mm]
		$^a$Department of Mathematics, Faculty of Basic Sciences, University of Qom, Qom, Iran.\\
		$^b$Shahab Danesh University, Qom, I. R. Iran.\\
		$^c$Department of Mathematics and Statistics, University of Maryland, Baltimore County, Baltimore, MD 21250, U.S.A.\\[2mm]
		{\tt E-mail: hassan.rezapour@gmail.com}\\
		{\tt E-mail: R.Nasiri@Shahabdanesh.ac.ir}\\
		{\tt E-mail: smousav1@umbc.edu}
	} \\[2mm]
\end{center}%
\vspace*{0.5cm}
\begin{quotation}
	\noindent
	{\footnotesize
		{\sc Abstract.}
Applications in chemistry motivated mathematicians to define different topological indices for different types of graphs. The Hyper-Zagreb index ($HM$) is an important tool as it integrates the first and the second Zagreb indices. In this paper, we characterize the trees and unicyclic graphs with the first four and first eight greatest $HM$-value, respectively.
}
\end{quotation}
\ \\
{\bf Keywords:} Hyper-Zagreb index, Vertex degree, Unicyclic graphs, Trees.\\

\n \textbf{2000 Mathematics subject classification: } 05C05, 05C07, 05C35.




\section{Introduction}
A nonnegative number can be assigned to a graph G to define an associated topological index if it is the same for every isomorphic graph of G, i.e., it is graph invariant. Topological indices are considered as appropriate tools to mathematically investigate and properly comprehend molecular structures and their properties such as complexity \cite{Gao_2017_Molecular_Descriptors, Gao_2016}. The first topological index was proposed by Wiener \cite{Wiener_1947} in order to examine chemical features of paraffin. Since trees turn out to have a special importance in various applications, authors in \cite{Dobrynin_2001} specifically study this index for these types of graphs. Moreover, In \cite{Nasiri2013}, the extremal unicyclic graphs with respect to Wiener index is studied. The Hyper-Wiener index for acyclic structures is due to Randic, where later \cite{Klein_1995} extends this notion so it can be applied for any connected graphs. An interested reader can explore some chemical applications of the Hyper-Wiener index in \cite{ Gutman_1997}. 

Zagreb index was first suggested by Gutman et al. \cite{Gutman_1972} in the 1970s, which absorbed attention of many scientists in different fields. A large amount of research has been done on this topic and the reader is encouraged to consult with \cite{Baby_2016, Braun_2005, Gutman_Das_2004, Kang_2016, Sonja_2003, Zhou_2004_Zagrebindices, Zhou_2005} for more useful information. A nice study on relations between the mentioned indices can be found in \cite{Zhou_2004}. 
\par
All graphs in this paper are assumed to be simple, finite and unidirected. The vertex and edge sets of a graph $G$ are shown by $V(G)$ and $E(G)$, respectively. Also, the number of vertices of $G$ is denoted by $n(G)$, which is called its order.
\par 
For a graph $G$, the \emph{Hyper-Zagreb index} of $G$ is defined as the following 
\begin{equation} \label{def: HZ-index}
HM(G)=\sum _{xy\in E(G)} \left( d_G(x)+d_G(y)\right)^2,
\end{equation}
where $d_G(x)$ is the degree of vertex $x$. 
For the edge $xy \in E(G)$, if consider $h_G(xy):=\left(d_G(x)+d_G(y)\right)^2$. Then, the above formulation can be equivalently written as 
\begin{equation*}
HM(G)=\sum _{xy\in E(G)} h_G(xy).
\end{equation*}
This invariant of graphs was initially presented by Shirdel et al. \cite{shirdel_2013} in 2013. They consider two simply connected graphs and compute this distance-based index for the resulted cartesian product, composition, join and disjunction graphs. Gao et al. \cite{Gao_2017} discuss acyclic, unicyclic, and bicyclic graphs and find sharp bounds for their Hyper-Zagreb index. The degree of vertices is the main part of some other graph invariants such as irregularity and total irregularity, see \cite{Fath-Tabar2013, NasiriEllahi, Nasiri-Gholami2017, Nasiri2014}.
To become more familiar with this topic, one should go through related literature including \cite{Basavanagoud_2017, FALAHATI-NEZHAD_2016, Siddiqui_2016, Kulli_2016}.
\section{Preliminaries and lemmas}
In this section, we first declare some basic and useful notations and definitions used in our work. Then, we propose some effective propositions which are important for achieving the goals of this article.
\par
\emph{Unicylcic graph} $G$ of order $n$ with circuit $C_m=x_1x_2\dots x_mx_1$ of length $m$ is denoted by $C_m^{u_1,u_2,\dots,u_k}\left(T_1,T_2,\dots, T_k\right)$ in which trees $T_i$'s for $i=1,2,\ldots, k$ are all nontrivial components of $G-E(C_m)$ and $u_i$ $\left(i=1,2,\dots,k\right)$ is the common vertex of $T_i$ and $C_m$. Specially, $G=C_n$ for
$k=0$. For convenience, we denote $C_m^{u_1,u_2,\ldots,u_k}\left(T_1,T_2,\ldots,T_k\right)$ by $C_m\left(T_1,T_2,\ldots,T_k\right)$, for any integer number $k\geq 1$. Let $n(T_i)=l_i+1$, $i=1,2,\ldots, k$, 
then $l=\sum_{i=1}^k l_i=n-m$. Also, if a tree $T_i$ is the star $S_{l_i+1}$ then we replace it by $l_i$, for example we denote $C_4\left(T_1,S_5,T_3,S_9\right)$ by $C_4\left(T_1,4,T_3,8\right)$. 
\par
\begin{figure}
	\begin{center}
		\begin{framed}
			\begin{tikzpicture}
			\def\r{.75}
			\def\l{2}
			\coordinate (x) at (0,0);
			\coordinate (y1) at (-2.5,-1);
			\coordinate (y2) at (-0.5,-1);
			\coordinate (yD) at (2.5,-1);
			\coordinate (c1) at (-2.5,-1-\r);
			\coordinate (c2) at (-0.5,-1-\r);
			\coordinate (cD) at (2.5,-1-\r);
			\draw[dashed] (c1) circle (\r );
			\draw[dashed] (c2) circle (\r );
			\draw[dashed] (cD) circle (\r );
			\filldraw (x) circle (1.5pt);
			\filldraw (y1) circle (1.5pt);
			\filldraw (y2) circle (1.5pt);
			\filldraw (yD) circle (1.5pt);
			\node at (c1){$T^*_1$};
			\node at (c2){$T^*_2$};
			\node at (cD){$T^*_\Delta$};
			\node [above] at (x){$x$};
			\node [anchor=north] at (y1){$y_1$};
			\node [anchor=north] at (y2){$y_2$};
			\node [anchor=north] at (yD){$y_\Delta$};
			\draw (x)--(y1);
			\draw (x)--(y2);
			\draw (x)--(yD);
			\draw[draw=none](c2)--node[sloped]{$\dots$}(cD);
			\end{tikzpicture}
		\end{framed}
	\end{center}
	\caption{Tree $T^x\left(T_1,T_2, \dots, T_{\Delta}\right)$.}\label{T*}
\end{figure}
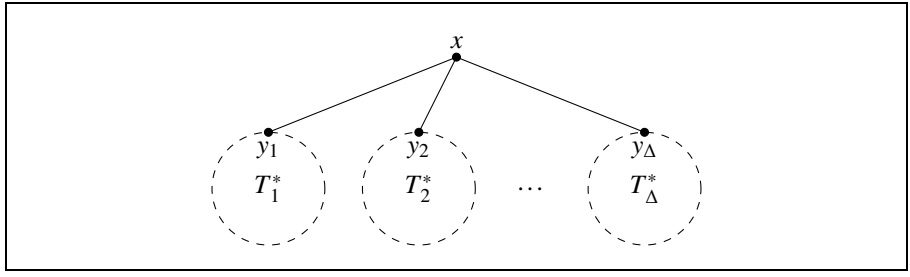

Let $T$ be a \emph{tree} with $n$ vertices $\left(n\ge 2\right)$ such that $x \in V(T)$ and $x$ has a maximum degree of vertices in graph $T$, i.e. $\Delta=d_T(x)=\max\left\{d_T(u), u \in V(T)\right\}$. $T$ is shown by $T^x\left(T_1,T_2, \dots, T_{\Delta}\right)$, where $T_i=T^*_i+\{y_ix\}$, $i=1,2, \dots, \Delta$, and $T^*_1, T^*_2, \dots, T^*_{\Delta}$
are trees with disjoint vertex sets and $n_1,n_2,\dots, n_{\Delta}$ are numbers of their vertices, respectively. 
Therefore, we have $|V(T_i)|=|V(T^*_i)|+1=n_i+1$, $i=1,2, \dots, \Delta$, and  $n=|V(T)|=\sum_{i=1}^{\Delta} n_i+1$ and $y_i\in V(T^*_i)$. Moreover, $E(T_i)=E(T^*_i)\cup \{y_ix\}$ and $V(T_i)=V(T^*_i)\cup \{x\}$ \ (see Figure \ref{T*}). 
\par
The \emph{coalescence} of G and H is denoted by $G(u)o H(v)$ and obtained by identifying
the vertex $u$ of $G$ with the vertex $v$ of $H$.
\begin{lemma} \label{lem: lem4}
	Assume that $z\in V(H)$ and $\{u,w\}\subseteq V(G)$ such that the following conditions hold:
	\begin{flalign*}
	\textit{(a)}&{}\,\, d_G(u)\le d_G(w),\\
	\textit{(b)}&{}  \sum_{x\in N_G(u)\setminus \{w\}}d_G(x) \le \sum_{x\in N_G(w)\setminus \{u\}}d_G(x).&&
	\end{flalign*}
		Moreover, let $G_1=G(u)o H(z)$ and $G_2=G(w)o H(z)$, where $G_1$ and $G_2$ are as shown in Figure \ref{Figure_2} . Then, $HM(G_2)\ge HM(G_1)$, with the equality if and only if  equality holds in both given conditions.
\end{lemma}
\begin{figure}[h] 
	\begin{center}
		\begin{framed}
			\begin{tikzpicture}
			\def\r{0.75}
			
			\begin{scope}[xshift=-3.5 cm]
			\coordinate (c1) at (0,0);
			\coordinate (c2) at (30:2*\r);
			\coordinate (u) at (30:\r);
			\coordinate (w) at  (-30:\r);
			\draw (c1) circle (\r );
			\draw (c2) circle (\r );
			\filldraw (u) circle (1.5pt);
			\filldraw (w) circle (1.5pt);
			\node at (c1){$G$};
			\node at (c2){$H$};
			\node [right] at (u){$u=z$};
			\node [right] at (w){$w$};
			\node at (1,-1.5*\r) {$G_1$};
			\end{scope}
			
			\begin{scope}[xshift=1.5 cm, yshift=0.75cm]
			\coordinate (c1) at (0,0);
			\coordinate (c2) at (-30:2*\r);
			\coordinate (u) at (30:\r);
			\coordinate (w) at  (-30:\r);
			\draw (c1) circle (\r );
			\draw (c2) circle (\r );
			\filldraw (u) circle (1.5pt);
			\filldraw (w) circle (1.5pt);
			\node at (c1){$G$};
			\node at (c2){$H$};
			\node [right] at (u){$u$};
			\node [right] at (w){$w=z$};
			\node at (1,-1.5*\r-1) {$G_2$};
			\end{scope}
			\end{tikzpicture}
		\end{framed}
	\end{center}
	\caption{The transformation of two graphs.}\label{Figure_2} 
\end{figure}
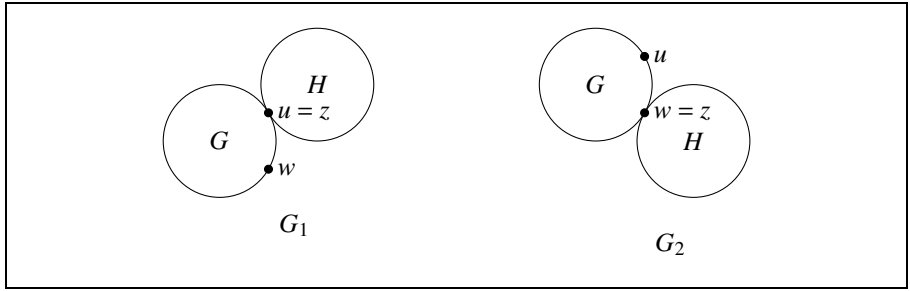
\begin{proof}
	Recall that 
	\begin{align*}
	\sum_{x \in N_G(w)\setminus \{u\}} h_{G_1}(xw)=&{} \sum_{x\in N_G(w) \setminus \{u\}}\left(d_G(w)+d_G(w)\right)^2,\\
	\sum_{x\in N_H(z)}h_{G_1}(xz)=&{}\sum_{x\in N_H(z)}\left(d_H(z)+d_G(u)+d_H(x)\right)^2
	\end{align*}
	and $HM_{G_1}(uw)=HM_{G_2}(uw)=\left(d_G(z)+d_G(w)+d_H(x)\right)^2.$
	In addition, one has
	\begin{align*}
	\sum_{x\in N_H(z)}h_{G_2}(zx)=&{}\sum_{x\in N_H(z)}\left(d_H(z)+d_G(w)+d_H(x)\right)^2,\\
	\sum_{x\in N_G(u)\setminus \{w\}}h_{G_2}(xu)=&{}\sum_{x\in N_G(u)\setminus \{w\}} \left(d_G(u)+d_G(x)\right)^2
	\end{align*}
	and 
	\begin{align*}
	\,\,\,\quad	\sum_{x\in N_G(w)\setminus \{u\}}h_{G_2}(xw)=\sum_{x\in N_G(w)\setminus \{x\}} \left(d_G(w)+d_G(x)+d_H(u)\right)^2.
	\end{align*}
	We consider two cases where either $uw \in E(G)$ or $uw \not \in E(G)$. First, suppose  that  $uw\in E(G)$. For $i=1$ and $2$, we have 
	\begin{align*}
	HM(G_i)=&{}\sum\limits_{\substack{xy \in E(G)\\ \ x,y\not \in \{u,v\}}} h_{G}(xy)+\sum _{x\in N_G(u)\setminus \{w\}}h_{G_i}(xu)+\sum_{x\in N_G(w)\setminus \{u\}}h_{G_i}(xw)\\
	+&{}HM_{G_i}(uw)+\sum_{x,y \not =z}h_H(xy) +\sum_{x\in N_H(z)}h_{G_i}(xz).
	\end{align*}
	On the other hand, 
	\begin{align*} 
	\sum_{N_G(u)\setminus \{w\}}h_{G_1}(xu)=\sum_{N_G(u)\setminus \{w\}}\left(d_G(u)+d_G(x)+d_H(x)\right)^2.
	\end{align*}
	Therefore,
	\begin{align*}
	HM(G_2)-HM(G_1)
	=&{}\sum_{x\in N_G(u)\setminus \{w\}}\left(\left(d_G(u)+d_G(x)\right)^2-\left(d_G(u)+d_G(x)+d_H(z)\right)^2\right)\\
	+&{}\sum_{x\in N_G(w)\setminus \{u\}} \left(\left(d_G(w)+d_G(x)+d_H(z)\right)^2-\left(d_G(w)+d_G(u)\right)^2\right)\\
	+&{}\sum_{x\in N_G(u)\setminus \{w\}}\left(\left(d_H(z)+d_G(w)+d_H(x)\right)^2-\left(d_H(z)+d_G(u)+d_H(w)\right)^2\right)
	\end{align*}
	this implies that
	\begin{align*}
	HM(G_2)-HM(G_1) \geq&{} 2d_H(z)\left(d_G(u)\left(d_G(w)-1\right)-d_G(u)\left(d_G(u)-1\right)\right)\\
	+&{}2d_H(z)\left(\sum_{x\in N_G(w)\setminus \{u\} } d_G(u)-\sum_{x\in N_G(w)\setminus \{w\}}d_G(x)
	\right)\\
	\ge&{} 0.
	\end{align*}
	Now, suppose that $uw\not \in E(G)$. Then, for $i=1$ and $2$, we have 
	\begin{align*}
	HM(G_i)=&{}\sum\limits_{\substack{xy\in E(G)\\ \ x,y \not \in \{u,w\}}} h_G(xy)+ \sum_{x\in N_G(u)} h_{G_i}(xu)+\sum_{x\in N_{G}(w)}h_{G_i}(xw)+\sum_{x,y \not =z}h_H(xy)\\
	+&{}\sum_{x\in N_{G_i}}h_{G_i}(xz).
	\end{align*}
	Also, in this case one has 
	$$\sum_{x\in N_G(w)\setminus \{u\}}d_G(x)=\sum_{x\in N_G(w)}d_G(x),
	\sum_{x\in N_G(u)\setminus \{w\}}d_G(x)=\sum_{x\in N_G(u)}d_G(x).$$
	Hence, a similar approach as the previous case can be used to prove the result. 
\end{proof}
\begin{lemma} \label{lem: lem1}
	Suppose $u$ and $v$ are vertices of graphs $G_1$ and $G_2$, respectively. Let $G$ be the graph obtained by joining $u\in V(G_1)$ to $v\in V(G_2)$ by an edge, and $G'$ be the graph obtained by identifying $u \in V(G_1)$ with $v\in V(G_2)$ and attaching a pendent vertex to the common vertex as shown in  Figure \ref{Figure_1}. Then if $d_G(u), d_{G'}(v)\ge 2$, we have $ HM(G)<HM(G')$. 
\end{lemma}
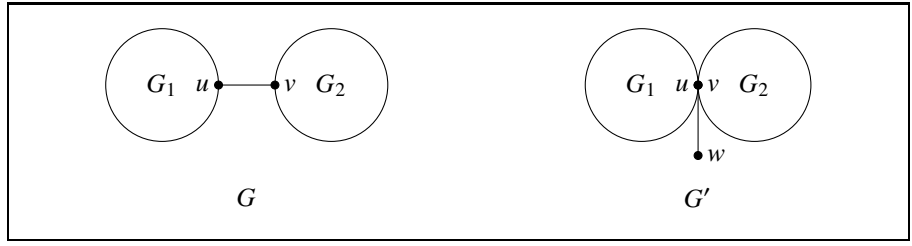
\begin{figure}[h]
	\begin{center}
		\begin{framed}
			\begin{tikzpicture}
			\usetikzlibrary{calc}
			\def\r{0.75}
			
			\begin{scope}[xshift=-4 cm]
			\coordinate (c1) at (-1.5*\r,0);
			\coordinate (c2) at (1.5*\r,0);
			\coordinate (u) at (-0.5*\r,0);
			\coordinate (v) at (0.5*\r,0);
			\draw (c1) circle (\r );
			\draw (c2) circle (\r );
			\draw (u)--(v);
			\filldraw (u) circle (1.5pt);
			\filldraw (v) circle (1.5pt);
			\node at (c1){$G_1$};
			\node at (c2){$G_2$};
			\node [left] at (u){$u$};
			\node [right] at (v){$v$};
			\node at (0,-2*\r) {$G$};
			\end{scope}
			
			\begin{scope}[xshift=2cm]
			\coordinate (c1) at (-\r,0);
			\coordinate (c2) at (\r,0);
			\coordinate (u) at (0,0);
			\coordinate (v) at (0,0);
			\coordinate (w) at (0,-1.25*\r);
			\draw (c1) circle (\r cm);
			\draw (c2) circle (\r );
			\draw (u)--(v);
			\filldraw (u) circle (1.5pt);
			\filldraw (v) circle (1.5pt);
			\filldraw (w) circle (1.5pt);
			\node at (c1){$G_1$};
			\node at (c2){$G_2$};
			\node [left] at (u){$u$};
			\node [right] at (v){$v$};
			\node [right] at (w){$w$};
			\draw (u)--(w);
			\node at (0,-2*\r) {$G'$};
			\end{scope}
			\end{tikzpicture}
		\end{framed}
	\end{center}
	\caption{ An illustration of graphs in Lemma \ref{lem: lem1}.}\label{Figure_1}
\end{figure}
\begin{proof}
	Assume that the graph $G'$ is obtained by identifying $u\in V(G_1)$ with $v\in V(G_2)$ and attaching a pendent vertex $w$ to the common vertex. Then,
\begin{align*}	
	HM(G)=&{} h_G(uv)+\sum_{x \in N_{G_1}(u)}h_G(ux)+\sum_{x\in N_{G_2}(v)} h_G(vx)
	+\sum\limits_{\substack{xy\in E(G_1)\\u \notin \{x,y\}}} h_{G_1}(xy)\\
	+&{}\sum_{xy\in E(G_2),v \notin \{x,y\}} h_{G_2}(xy)
\end{align*}
	and
\begin{align*}
	HM(G')=&{} h_{G'}(uw)+\sum_{x \in N_{G_1}(u)}h_{G'}(ux)
	+
	\sum_{x\in N_{G_2}(v)} h_{G'}(vx)+\sum\limits_{\substack{xy\in E(G_1)\\u \notin \{x,y\}}} h_{G_1}(xy)\\
	+&{}
	\sum\limits_{\substack{xy\in E(G_2)\\v \notin \{x,y\}}} h_{G_2}(xy).
\end{align*}
	Since $d_G(u)=d_{G_1}(u)+1, \ d_G(v)=d_{G_1}(v)+1, \ d_{G'}(w)=1$ and $d_{G'}(u)=d_{G'}(v)=d_{G_1}(u)+d_{G_2}(v)+1$ we have
\begin{align*}
	\sum_{x \in N_{G_1}(u)} h_G(ux)< \sum_{x \in N_{G_1}(u)} h_{G'}(ux),\quad \quad
	\sum_{x \in N_{G_2}(v)} h_G(vx)< \sum_{x \in N_{G_2}(v)} h_{G'}(vx) 
\end{align*}
	and $h_G(uv)=h_{G'}(uw)=\left(d_{G_1}(u)+d_{G_2}(v)+2\right)^2$. Hence, 
\begin{align*}	
	HM(G')-HM(G)=&{}
	\sum_{x \in N_{G_1}(u)} h_{G'}(ux)
	-
	\sum_{x \in N_{G_1}(u)} h_{G}(ux) \\
	+&{}
	\sum_{x \in N_{G_2}(v)} h_{G'}(vx)
	-
	\sum_{x \in N_{G_2}(v)} h_{G}(vx)\\
	>&{}0.
\end{align*}
\end{proof}


\begin{corollary} \label{cor: cor0}
	Let $T$ be a tree with $n$ vertices. Then, $HM(T)\le HM\left(S_n\right)$, with the equality if and only if $T \cong S_n$.
\end{corollary}

\begin{corollary} \label{cor: cor00}
	Let $G=C_m\left(T_1,T_2,\dots, T_k\right)$ be a unicyclic graph and $n(T_i)=l_i+1$. Then, 
	$HM(G)\leq HM\left(C_m\left(l_1,l_2,\dots, l_k\right)\right)$, 
	with the equality if and only if $T_i \cong S_{l_i+1}$, $i=1,2,\dots, k$.
\end{corollary}
\begin{lemma} \label{lem: lem5}
	Let $G_1=C_m\left(l_1,l_2,\dots, l_k\right)$ be a unicyclic graph and $y_1u_i, 	u_iu_{i+1} \in E(C_m)$ such that $d_{G_1}(y_1),d_{G_1}(u_1) \le d_{G_1}(u_{i+1})$, then for  $G_2=C_m\left(l_1,\dots, l_{i-1},l_{i+1}+l_i,l_{i+2},\dots,l_k\right)$ one has that  $HM(G_1)< HM(G_2)$.
\end{lemma}

\begin{proof}
	Let $G=C_m\left(l_1,\dots, l_{i-1},l_{i+1},l_{i+2},\dots,l_k\right)$, then $2=d_G(u_i)<3\le d_G(u_{i+1})$; meaning that the condition \textit{(a)} in Lemma \ref{lem: lem4} holds. Hence, we now show that the second condition in this Lemma is also satisfied. Suppose that $y_2u_{i+1}\in E(C_m)$. By a simple calculation one can check that
	\begin{align*}
	\sum_{x \in N_G(u_i)\setminus \{u_{i+1} \}} d_G(x)
	=&{}d_G(y_1),\\
	\sum_{x \in N_G(u_{i+1})\setminus \{u_{i} \}} d_G(x)
	=&{}\sum\limits_{\substack{x \in N_G(u_{i+1})\setminus \{u_{i} \}\\x\in V(C_m)}} d_G(x)+\sum\limits_{\substack{x \in N_G(u_{i+1})\setminus \{u_{i} \}\\x \notin V(C_m)}} d_G(x)\\
	=&{}d_G(y_2)+\sum\limits_{\substack{x \in N_G(u_{i+1})\setminus \{u_{i} \}\\x \notin V(C_m)}} 1\\
	=&{}d_G(y_2)+d_G(u_{i+1})-2.
	\end{align*}
Moreover, $d_{G_1}(u_{i+1})=d_G(u_{i+1})$ and $d_G(y_1)=d_{G_1}(y_1)$. On the other hand, since $y_2 \in V(C_m)$ then $d_G(y_2)\ge 2$; implying that $d_{G}(y_2)-2\ge 0$. So, we have 
	\begin{align*}
\sum_{x \in N_G(u_i)\setminus \{u_{i+1}\} } d_G(x)=d_G(y_1)=d_{G_1}(y_1)
\le d_{G_1}(u_{i+1})= d_G(u_{i+1})\le \sum_{x \in N_G(u_{i+1})\setminus \{u_i\}} d_G(x).
	\end{align*}
	Therefore, the condition \textit{(b)} of Lemma \ref{lem: lem4} holds, which completes the proof.
\end{proof}

%
\begin{lemma} \label{lem: lem6}
	Let $G=C_m^{u_1,u_2, \dots, u_k}\left(l_1,l_2,\dots, l_k\right)$ be a unicyclic graph and $k>1$. Then if $u_iu_{i+1}\in E(C_m)$, $i=1,2,\dots,k-1$, then $HM(G)< HM\left(C_m\left(n-m\right)\right)$. Otherwise, there exist positive integers $l'_1, l'_2, \dots, l'_r \ (r\le k)$, such that $HM(G)<HM(G')<HM(G'')$, where $G'=C_m^{v_1,v_2, \dots, v_r}\left(l'_1, l'_2, \dots, l'_r\right)$, $G''=C_m^{v_2,v_3, \dots, v_r}\left(l'_1+l'_2, l'_3, \dots, l'_r\right)$, $d_{G'}(v_i, v_j) \ge 2$ \ for \ $1\le i < j \le r$ \ and \ $\{v_1, v_2, \dots, v_r\} \subseteq \{u_1, u_2, \dots, u_k\}$.
\end{lemma}

\begin{proof}
$HM(G)<HM(G')$ is straightforward in light of Lemma \ref{lem: lem5}. Now, by considering $u=v_1, w=v_2$, $H=S_{l'_1+1}$ and $G=C_m^{v_2,v_3, \dots, v_r}\left(l'_2, l'_3, \dots, l'_r\right)$ and using Lemma \ref{lem: lem4} we can conclude that $HM(G')<HM(G'')$, as desired.
\end{proof}

\begin{lemma} \label{lem: R}
Let $G=C_m\left(l_1,l_2,\dots, l_k\right)$ be a unicyclic graph of order $n$. Then $HM(G)\le HM\left(C_m\left(n-m\right)\right)$, with equality if and only if $k=1$.	
\end{lemma}
\begin{proof}
The proof is obtained by applying Lemmas \ref{lem: lem5} and \ref{lem: lem6}.
\end{proof}
\begin{corollary} \label{cor: cor1}
Let $G=C_m\left(T_1,T_2,\dots, T_k\right)$ be a unicyclic graph and $n(T_i)=l_i+1$, $i=1,2,\dots,k$. Then, 
$$HM(G)\leq HM\left(C_m\left(l_1,l_2,\dots, l_k\right)\right)\le HM\left(C_m\left(n-m\right)\right),$$ 
with left equality if and only if $T_i \cong S_{l_i+1}$, $i=1,2,\dots, k$, and right equality if and only if $k=1$.
\end{corollary}

\begin{lemma} \label{lem: newlem}
	Let $G_1=C_m\left(n-m\right)$ and $G_2=C_{m-1}\left(n-m+1\right)$, $m\geq 4$, be unicyclic graphs of order $n$. Then, $HM(G_1)<HM(G_2)$.
\end{lemma}
\begin{proof}
By a simple calculation we have
\begin{align*}
HM(G_1)
=&{}4\left(m-2\right)+2\left(n-m+4\right)^2+(n-m)\left(n-m+3\right)^2\\
<&{}4\left(m-3\right)+2\left(n-m+5\right)^2+(n-m+1)\left(n-m+4\right)^2\\
=&{}HM(G_2).
\end{align*}
As desired.
\end{proof}
\begin{lemma} \label{lem: lem3}
	Let $G=C_m\left(T_1,T_2,\dots, T_k\right)$ be a unicyclic graph. Then, 
	\begin{align*}
	HM(G)=\sum_{i=1}^k\sum_{xy \in E(T_i)}h_{G}(xy)+\sum_{xy \in E(C_m)}h_{G}(xy).
	\end{align*}
\end{lemma}
\begin{proof}
	The proof is trivial by the Hyper-Zagreb index definition \eqref{def: HZ-index}.
\end{proof}



\section{Main Results}

\begin{figure}[h]
	\begin{center}
		\begin{framed}
			\resizebox{\textwidth}{!}{
				\begin{tikzpicture}
				\def\l{1}
				\def\r{2}
				\newcommand{\txtsize}{\scriptsize}
				\begin{scope}
				\coordinate (n) at (0:0);
				\coordinate (1) at (0:\l);
				\coordinate(2) at (60:\l);
				\coordinate(3) at (120:\l);
				\begin{scope}[shift={(120:\l)}]
				\coordinate(4)at (60:\l);
				\end{scope}
				\draw[draw=white, line width=0pt](1)--node[sloped]{$\cdots$} (2) ;
				\draw [decorate,decoration=brace] ([shift= (60:\l)]70:0.2) -- node [sloped, black,midway,yshift=8pt] {\txtsize$n-3$} ([shift= (0:\l)]-10:0.2);
				\draw (1)--(n)--(2);
				\draw (4)--(3)--(n);
				\foreach \i in {n,1,2,...,4}{
					\filldraw (\i) circle (\r pt);
				}
				\node at(0.5*\l,-0.8*\l){\txtsize $T^1_n$};
				\end{scope}
				\begin{scope}[shift={(0:3.8*\l)}]
				\coordinate (n) at (0:0);
				\coordinate (1) at (0:\l);
				\coordinate(2) at (60:\l);
				\coordinate(3) at (120:\l);
				\begin{scope}[shift={(120:\l)}]
				\coordinate(4)at (60:\l);
				\coordinate (5) at (120:\l);
				\end{scope}
				\draw[draw=white, line width=0pt](1)--node[sloped]{$\cdots$} (2) ;
				\draw [decorate,decoration=brace] ([shift= (60:\l)]70:0.2) -- node [sloped, black,midway,yshift=8pt] {\txtsize$n-4$} ([shift= (0:\l)]-10:0.2);
				\draw (1)--(n)--(2);
				\draw (4)--(3)--(n);
				\draw (5)--(3);
				\foreach \i in {n,1,2,...,5}{
					\filldraw (\i) circle (\r pt);
				}
				\node at(0.5*\l,-0.8*\l){\txtsize $T^2_n$};
				\end{scope}
				\begin{scope}[shift={(0:8.1*\l)}]
				\coordinate (n) at (0:0);
				\coordinate (1) at (0:\l);
				\coordinate(2) at (60:\l);
				\coordinate(3) at (120:\l);
				\begin{scope}[shift={(120:\l)}]
				\coordinate(4)at (120:\l);
				\end{scope}
				\coordinate (5) at (180:\l);
				\begin{scope}[shift={(180:\l)}]
				\coordinate(6)at (120:\l);
				\end{scope}
				\draw[draw=white, line width=0pt](1)--node[sloped]{$\cdots$} (2) ;
				\draw [decorate,decoration=brace] ([shift= (60:\l)]70:0.2) -- node [sloped, black,midway,yshift=8pt] {\txtsize$n-5$} ([shift= (0:\l)]-10:0.2);
				\draw (1)--(n)--(2);
				\draw (4)--(3)--(n);
				\draw (6)--(5)--(n);
				\foreach \i in {n,1,2,...,6}{
					\filldraw (\i) circle (\r pt);
				}
				\node at(0,-0.8*\l){\txtsize $T^3_n$};
				\end{scope}
				\begin{scope}[shift={(0:12.1*\l)}]
				\coordinate (n) at (0:0);
				\coordinate (1) at (0:\l);
				\coordinate(2) at (60:\l);
				\coordinate(3) at (120:\l);
				\begin{scope}[shift={(120:\l)}]
				\coordinate (5) at (60:\l);
				\coordinate(4)at (120:\l);
				\coordinate(6)at (-120:\l);
				\end{scope}
				\draw[draw=white, line width=0pt](1)--node[sloped]{$\cdots$} (2) ;
				\draw [decorate,decoration=brace] ([shift= (60:\l)]70:0.2) -- node [sloped, black,midway,yshift=8pt] {\txtsize$n-6$} ([shift= (0:\l)]-10:0.2);
				\draw (1)--(n)--(2);
				\draw (4)--(3)--(n);
				\draw (5)--(3)--(6);
				\foreach \i in {n,1,2,...,6}{
					\filldraw (\i) circle (\r pt);
				}
				\node at(0,-0.8*\l){\txtsize $T^4_n$};
				\end{scope}
				\end{tikzpicture}
			}
		\end{framed}
	\end{center}		
	\caption{Some trees with large Hyper-Zagreb values.}\label{Figure_4}
\end{figure}
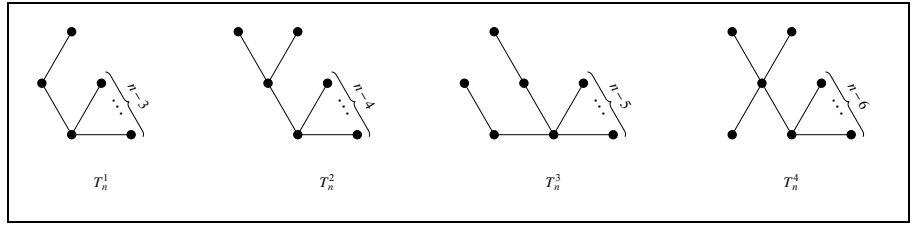

\begin{figure}[h]
	\begin{center}
		\begin{framed}
			\resizebox{\textwidth}{!}{
				\begin{tikzpicture}
				\def\l{1}
				\def\r{2}
				\newcommand{\txtsize}{\scriptsize}
				\begin{scope}[shift={(180:0.5*\l)}]
				
				\begin{scope}[shift={(0,0)}]
				\coordinate (n-2) at (0,0);
				\coordinate (n-1) at (-120:\l);
				\coordinate (n) at (-60:\l);
				\coordinate (1) at (150:\l);
				\coordinate (n-5) at (210:\l);
				\draw[draw=white, line width=0pt](1)--node[sloped]{$\cdots$} (n-5) ;
				\draw [decorate,decoration=brace] ([shift= (210:\l)]200:0.2) -- node [sloped, black,midway,yshift=8pt] {\txtsize$n-3$} ([shift=(150:\l)]160:0.2);
				\foreach \i in {1,n-5,n-1,n}{
					\draw (\i)--(n-2);
					\filldraw (\i) circle (\r pt);
				}
				\draw (n)--(n-1);
				\foreach \i in {n-2}{
					\filldraw (\i) circle (\r pt);
				}
				\node at(0,-1.7*\l){\txtsize $C_3\left(n-3\right)$};
				\end{scope}
				
				\begin{scope}[shift={(4,0)}]
				\coordinate (n-2) at (0,0);
				\coordinate (n-1) at (-120:\l);
				\coordinate (n) at (-60:\l);
				\coordinate (1) at (150:\l);
				\coordinate (n-5) at (210:\l);
				\coordinate (n-4) at (90:\l);
				\coordinate (n-3) at ([shift=(90:\l)]90:\l);
				\draw[draw=white, line width=0pt](1)--node[sloped]{$\cdots$} (n-5) ;
				\draw [decorate,decoration=brace] ([shift= (210:\l)]200:0.2) -- node [sloped, black,midway,yshift=8pt] {\txtsize$n-5$} ([shift=(150:\l)]160:0.2);
				\foreach \i in {1,n-5,n-4,n-1,n}{
					\draw (\i)--(n-2);
					\filldraw (\i) circle (\r pt);
				}
				\draw (n-3)--(n-4);
				\draw (n)--(n-1);
				\foreach \i in {n-2,n-3}{
					\filldraw (\i) circle (\r pt);
				}
				\node at(0,-1.7*\l){\txtsize $C_3\left(T^1_{n-2}\right)$};
				\end{scope}
				\begin{scope}[shift={(8,0)}]
				\coordinate (n-2) at (0,0);
				\coordinate (n-1) at (-120:\l);
				\coordinate (n) at (-60:\l);
				\coordinate (1) at (150:\l);
				\coordinate (n-6) at (210:\l);
				\coordinate (n-5) at (90:\l);
				\coordinate (n-3) at (-0.6,2*\l);
				\coordinate (n-4) at (0.6,2*\l);
				\draw (n-3)--(n-5)--(n-4);
				\draw[draw=white, line width=0pt](1)--node[sloped]{$\cdots$} (n-6) ;
				\draw [decorate,decoration=brace] ([shift= (210:\l)]200:0.2) -- node [sloped, black,midway,yshift=8pt] {\txtsize$n-6$} ([shift=(150:\l)]160:0.2);
				\foreach \i in {1,n-5,n-6,n-1,n}{
					\draw (\i)--(n-2);
					\filldraw (\i) circle (\r pt);
				}
				\draw (n)--(n-1);
				\foreach \i in {n-2,n-3,n-4}{
					\filldraw (\i) circle (\r pt);
				}
				\node at(0,-1.7*\l){\txtsize $C_3\left(T^2_{n-2}\right)$};
				\end{scope}
				\begin{scope}[shift={(12,0)}]
				\coordinate (n-2) at (0,0);
				\coordinate (n-1) at (-120:\l);
				\coordinate (n) at (-60:\l);
				\coordinate (1) at (150:\l);
				\coordinate (n-7) at (210:\l);
				\coordinate (n-4) at (0.6,\l);
				\coordinate (n-6) at (-0.6,\l);
				\coordinate (n-5) at (-0.6,2*\l);
				\coordinate (n-3) at (0.6,2*\l);
				\draw (n-6)--(n-5);
				\draw (n-3)--(n-4);
				\draw[draw=white, line width=0pt](1)--node[sloped]{$\cdots$} (n-7) ;
				\draw [decorate,decoration=brace] ([shift= (210:\l)]200:0.2) -- node [sloped, black,midway,yshift=8pt] {\txtsize$n-7$} ([shift=(150:\l)]160:0.2);
				\foreach \i in {1,n-4,n-6,n-1,n,n-7}{
					\draw (\i)--(n-2);
					\filldraw (\i) circle (\r pt);
				}
				\draw (n)--(n-1);
				\foreach \i in {n-2,n-3,n-4,n-5}{
					\filldraw (\i) circle (\r pt);
				}
				\node at(0,-1.7*\l){\txtsize $C_3\left(T^3_{n-2}\right)$};
				\end{scope}
				\end{scope}
				\begin{scope}[shift={(0,-4.5*\l)}]
				\begin{scope}[shift={(180:0.5*\l)}]
				\coordinate (n-2) at (0:0);
				\coordinate (n-1) at (0:\l);
				\coordinate (n) at (60:\l);
				\coordinate (1) at (120:\l);
				\coordinate (n-6) at (180:\l);
				\begin{scope}[shift={(60:\l)}]
				\coordinate (n-5) at (120:\l);
				\end{scope}
				\draw[draw=white, line width=0pt](1)--node[sloped]{$\cdots$} (n-6) ;
				\draw [decorate,decoration=brace] ([shift= (180:\l)]190:0.2) -- node [sloped, black,midway,yshift=8pt] {\txtsize$n-4$} ([shift= (120:\l)]110:0.2);
				\foreach \i in {1,n-1,n,n-6}{
					\draw (\i)--(n-2);
					\filldraw (\i) circle (\r pt);
				}
				\draw (n-5)--(n);
				\draw (n)--(n-1);
				\foreach \i in {n-2,n-5}{
					\filldraw (\i) circle (\r pt);
				}
				\node at(0,-0.8*\l){\txtsize $C_3\left(1,n-4\right)$};
				\end{scope}
				\begin{scope}[shift={(3.8,0)}]
				\coordinate (n-2) at (0:0);
				\coordinate (n-1) at (0:\l);
				\coordinate (n) at (60:\l);
				\coordinate (1) at (120:\l);
				\coordinate (n-5) at (180:\l);
				\begin{scope}[shift={(60:\l)}]
				\coordinate (n-4) at (60:\l);
				\begin{scope}[shift={(60:\l)}]
				\coordinate (n-3) at (180:\l);
				\end{scope}
				\end{scope}
				\draw[draw=white, line width=0pt](1)--node[sloped]{$\cdots$} (n-5) ;
				\draw [decorate,decoration=brace] ([shift= (180:\l)]190:0.2) -- node [sloped, black,midway,yshift=8pt] {\txtsize$n-5$} ([shift= (120:\l)]110:0.2);
				\foreach \i in {1,n-5,n-1,n}{
					\draw (\i)--(n-2);
					\filldraw (\i) circle (\r pt);
				}
				\draw (n-4)--(n)--(n-3);
				\draw (n)--(n-1);
				\foreach \i in {n-2,n-3,n-4}{
					\filldraw (\i) circle (\r pt);
				}
				\node at(0,-0.8*\l){\txtsize $C_3\left(2,n-5\right)$};
				\end{scope}
				\begin{scope}[shift={(7.6,0)}]
				\coordinate (n-2) at (0:0);
				\coordinate (n-1) at (0:\l);
				\coordinate (1) at (120:\l);
				\coordinate (n-6) at (180:\l);
				\begin{scope}[shift={(180:0.5*\l)}]
				\coordinate (n) at (60:\l);
				\coordinate (n-3) at ([shift=(0:\l)]60:\l);
				\end{scope}
				\draw[draw=white, line width=0pt](1)--node[sloped]{$\cdots$} (n-6) ;
				\draw [decorate,decoration=brace] ([shift= (180:\l)]190:0.2) -- node [sloped, black,midway,yshift=8pt] {\txtsize$n-4$} ([shift= (120:\l)]110:0.2);
				\foreach \i in {1,n-1,n,n-6}{
					\draw (\i)--(n-2);
					\filldraw (\i) circle (\r pt);
				}
				\draw (n)--(n-3)--(n-1);
				\foreach \i in {n-2,n-3}{
					\filldraw (\i) circle (\r pt);
				}
				\node at(0,-0.8*\l){\txtsize $C_4\left(n-4\right)$};
				\end{scope}
				\begin{scope}[shift={(11.4,0)}]
				\coordinate (n-2) at (0:0);
				\coordinate (n-1) at (0:\l);
				\coordinate (n) at (60:\l);
				\coordinate (1) at (120:\l);
				\coordinate (n-6) at (180:\l);
				\begin{scope}[shift={(60:\l)}]
				\coordinate (n-5) at (120:\l);
				\end{scope}
				\coordinate (n-3) at ([shift=(0:\l)]60:\l);
				\draw[draw=white, line width=0pt](1)--node[sloped]{$\cdots$} (n-6) ;
				\draw [decorate,decoration=brace] ([shift= (180:\l)]190:0.2) -- node [sloped, black,midway,yshift=8pt] {\txtsize$n-5$} ([shift= (120:\l)]110:0.2);
				\foreach \i in {1,n-1,n,n-6}{
					\draw (\i)--(n-2);
					\filldraw (\i) circle (\r pt);
				}
				\draw (n)--(n-5);
				\draw (n)--(n-1)--(n-3);
				\foreach \i in {n-2,n-3,n-5}{
					\filldraw (\i) circle (\r pt);
				}
				\node at(0,-0.8*\l){\txtsize $C_3\left(1,1,n-5\right)$};
				\end{scope}
				\end{scope}
				\end{tikzpicture}
			}
		\end{framed}
	\end{center}
	\caption{The unicyclic graphs with the first eight greatest Hyper-Zagreb.}\label{Uni}
\end{figure}
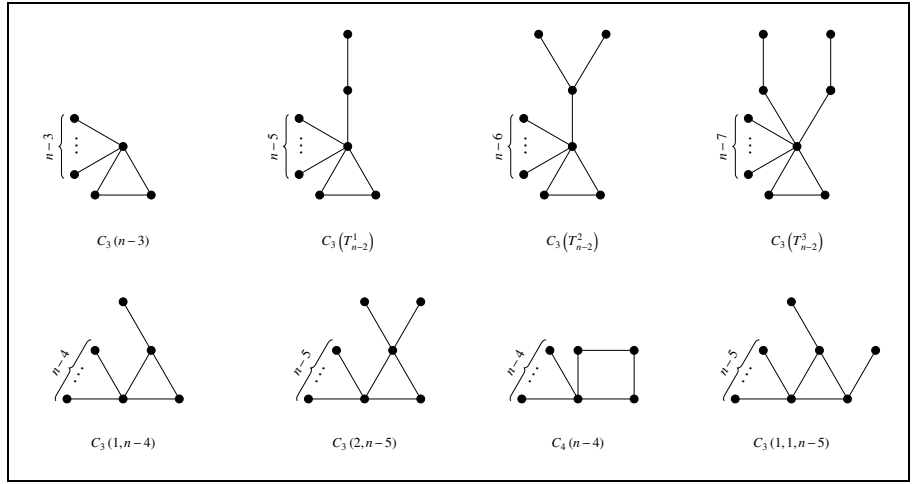

In this section, we characterize the trees and unicyclic graphs with the first four and first eight greatest HM-value, respectively.
\begin{theorem} \label{thm: thm1}
	Let $T$ be a tree with $n$ vertices. If $T\not \cong S_n, T^1_n$ or $T^n_2$, then 
	$$HM(T)\le HM\left(T^3_n\right)< HM\left(T^2_n\right)< HM\left(T^1_n\right)< HM\left(S_n\right),$$
	 with the equality if and if $T\cong T^n_3$, where $T^1_n,T^2_n$ and $T^3_n$ are given as in Figure \ref{Figure_4}.
\end{theorem}
\addtolength{\tabcolsep}{30pt}    
\begin{table}[h]
\begin{center}
	\begin{tabular}{ ||l|l||} 
		\multicolumn{2}{c} {}        \\  \hline
		Graph    &     HM-value   \\ \hline \hline
		$S_n$   & $n^3-n^2$          \\  \hline
		$T^n_1$ & $n^3-4n^2+7n+6$    \\ \hline 
		$T^n_2$ & $n^3-7n^2+20n+16$  \\  \hline
		$T^n_3$ & $n^3-7n^2+20n$     \\  \hline 
		\end{tabular}
	\caption{Trees with large Hyper-Zagreb values.}\label{tab: tab1}
\end{center}
\end{table}
\addtolength{\tabcolsep}{-30pt}
\begin{proof}
	Using Table \ref{tab: tab1}, we have $HM\left(S_n\right)> HM\left(T^1_n\right)> HM\left(T^2_n\right)> HM\left(T^3_n\right)$. Hence, we need to prove  that $HM(T)<HM\left(T^3_n\right)$ when $T\not \cong T^3_n$. Let $T=T^x\left(T_1, T_2,\dots, T_{\Delta}\right)$, where $\Delta= d_T(x)$. By Corollary \ref{cor: cor0}, we have $HM\left(T_i\right)\le HM\left(S_{n_i}\right)$,  $i=1,2,\dots,\Delta$. Moreover, let $T'=T^x\left(T'_1, T'_2, \dots. T'_{\Delta}\right)$, where $T'^*_i=S_{n_i}$, $i=1,2,\dots,\Delta$, then we have $HM(T)\le HM\left(T'\right)$. To complete the proof, we consider three different cases as follows:
\par
\textit{Case 1}: 
	assume that $d_T(y_i)=1$, $i=1,2,\dots, \Delta$, then $T=S_n$. This is a contradiction to the assumption. 
\par
\textit{Case 2}: 
	assume that there exists  $y_t$ for $t=1, 2, \dots, \Delta $ such that $d_T(y_t) \ge 2$ and $d_T(y_i)=1$  for $\ i=1,2,\ldots, \Delta$ and $i\not = t$. In this case, there are three subcases that can happen:
	\begin{enumerate}
\item[\textit{(i)}] If $\left|V(T^*_t)\right|=2$, then $T\cong T^n_1$. This is clearly a contradiction.
\item[\textit{(ii)}] If $\left|V(T^*_t)\right|=3$, then we must consider that $d_T(y_t)=2$ or $3$. The case $d_T(y_t)=3$ implies that $T\cong T^2_n$, which is a contradiction. If $d_T(y_t)=2$, then 
	\begin{align*}
	HM(T)=&{}\left(n-4\right)\left(n-2\right)^2+\left(n-1\right)^2+16+9\\
	=&{}n^3-7n^2+18n+10\\
	<&{}n^3-7n^2+20n\\
	=&{}HM\left(T^3_n\right).
	\end{align*}
\item[\textit{(iii)}] If $\left|V(T^*_t)\right|\ge 4 $, then $T=T^x\left(\overbrace{ \ S_2,\dots, S_2}^{t-1 \ times}, \  T_t, \overbrace{\ S_2, \dots, S_2}^{\Delta-t \ times }\right)$. By Corollary \ref{cor: cor0}, the Hyper-Zagreb index for $T$ is maximum when $T^*_t=S_{n_t}$. On the other hand, it follows from Lemma \ref{lem: lem4} that if $n_t=4$ than $T$ has maximum $HM$-value, i.e. in this case $T$ has maximum $HM$-value when $T\cong T^4_n$ (see Figure \ref{Figure_4}).
Hence, applying Lemma \ref{lem: lem4}, it is clear that $HM(T)\le HM\left(T^4_n\right)<HM\left(T^3_n\right)$.
	\end{enumerate}
\par
\textit{Case 3}: Suppose that there exist $1\le s,t\le \Delta$ such that $d_T(y_s), d_T(y_t)\ge 2$. Similar to previous, applying Corollary \ref{cor: cor0} and Lemma \ref{lem: lem4}, it can concluded that $HM(T)<HM\left(T^3_n\right)$.
\end{proof}
\begin{theorem} \label{thm: thm2}
	Let $G$ be a unicylcic graph of order $n\ge 15$. If 
$G\not \cong C_3\left(n-3\right),\\ \ C_3\left(1,n-4\right), \ C_3\left(T^1_{n-2}\right), \ C_4\left(n-2\right), \ C_3\left(2,n-5\right), \ C_3\left(1,1,n-5\right) \ and \ C_3\left(T^2_{n-2}\right)$.
	Then, 
	\begin{align*}
	HM(G)\le&{}
	HM\left(C_3\left(T^3_{n-2}\right)\right)<HM\left(C_3\left(T^2_{n-2}\right)\right)<HM\left(C_3\left(1,1,n-5\right)\right)\\
	<&{}HM\left(C_3\left(2,n-5\right)\right)<HM\left(C_4\left(n-4\right)\right) <HM\left(C_3\left(T^1_{n-2}\right)\right)\\
	<&{}H M\left(C_3\left(1,n-4\right)\right)< HM\left(C_3\left(n-3\right)\right),
	\end{align*}
	with the equality if and only if $G\cong C_3\left(T^3_{n-2}\right)$ or $G\cong C_3\left(P_3,10\right)$ for $n= 15$ (see Figure \ref{Uni}).
\end{theorem}
\addtolength{\tabcolsep}{30pt}    
\begin{table}[h]
\begin{center}
	\begin{tabular}{ ||l|l||} 
		\multicolumn{2}{c} {}  \\  \hline
		Graph    & HM-value   \\ \hline \hline
		$C_3\left(n-3\right)$         & $n^3-n^2+4n+18$          \\  \hline
		$C_3\left(1,n-4\right)$ & $n^3-4n^2+11n+38$    \\ \hline 
		$C_3\left(T^1_{n-2}\right)$       & $n^3-4n^2+11n+20$  \\  \hline
		$C_4\left(n-4\right)$         & $n^3-4n^2+9n+28$     \\  \hline
		$C_3\left(2,n-5\right)$     & $n^3-7n^2+24n+68$     \\  \hline
		$C_3\left(1,1,n-5\right)$ & $n^3-7n^2+24n+48$     \\  \hline
		$C_3\left(T^2_{n-2}\right)$                 & $n^3-7n^2+24n+26$     \\  \hline
		$C_3\left(T^3_{n-2}\right)$                  & $n^3-7n^2+24n+10$     \\  \hline
	\end{tabular}
	\caption{Unicyclic graphs with large Hyper-Zagreb values.}\label{tab: tab2}
\end{center}
\end{table}
\addtolength{\tabcolsep}{-30pt} 
\begin{proof}
	Assume that $G=C_m\left(T_1,T_2,\dots, T_k\right)$ be a unicyclic graph and $n(T_i)=l_i+1$, $i=1,2,\dots, k$.
	The given Table \ref{tab: tab2} provides the  Hyper-Zagreb index of some graphs by which the result is trivial. It is enough to discuss about the equality case. If $G\cong C_3\left(T^3_{n-2}\right)$, then $HM(G)=HM\left(C_3\left(T^3_{n-2}\right)\right)$. Also, if $G\cong C_3\left(P_3,10\right)$ for $n=15$, then $HM(G)=2170=HM\left(C_3\left(T^3_{13}\right)\right)$. We now prove that $HM(G)<HM\left(C_3\left(T^3_{n-2}\right)\right)$, where 
	$G\not \cong C_3\left(n-3\right), \ C_3\left(1,n-4\right), \ C_3\left(T^1_{n-2}\right), \ C_4\left(n-2\right), \ C_3\left(2,n-5\right), \ C_3\left(1,1,n-5\right)$ and $C_3\left(T^2_{n-2}\right)$.
	We examine three cases of $m=3, 4$ and $5$ for $G=C_m\left(T_1,T_2,\dots, T_k\right)$ as follows:
	\par
	\textit{Case 1}: $m=3$. We need to discuss three subcases that $k=1,2$ and $3$.
	\begin{enumerate}
		\item[\textit{(i)}] $k=1$, then $G=C_3\left(T_1\right)$. By assumption, we know that $T_1\not \cong S_{n-2}, T^1_{n-2}, T^2_{n-2}$ and $ T^3_{n-2}$. So, Theorem \ref{thm: thm1} implies that $HM(T_1)<HM(T^3_{n-2})$. Lemma \ref{lem: lem3} now guarantees that $HM(G)<HM\left(C_3\left(T^3_{n-2}\right)\right).$ 
		\item[\textit{(ii)}] $k=2$, then $G=C_3\left(T_1,T_2\right)$. By assumption, $G\not \cong C_3\left(1,n-4\right)$ and $C_3\left(2,n-5\right)$. By Corollaries \ref{cor: cor0}, \ref{cor: cor00} and Lemmas \ref{lem: lem4}, \ref{lem: lem3}, the maximum value of $HM(G)$ happens when  $G\cong C_3\left(3,n-6\right)$ or $C_3\left(P_3,n-5\right)$. The first case yields that
		\begin{align*}
         HM(G)\le HM\left(C_3\left(3,n-6\right)\right)=n^3-10n^2+43n+108.
		\end{align*}
	Hence, we have (for $n\ge 15$) that
		\begin{align*}
			HM\left(C_3\left(T^3_{n-3}\right)\right)-HM(G)
			\geq&{}HM\left(C_3\left(T^3_{n-3}\right)\right)-HM\left(C_3\left(3,n-6\right)\right)\\
		=&{}\left(n^3-7n^2+24n+10\right)\\
		-&{}\left(n^3-10n^2+43n+108\right)\\
		=&{}3n^2-19n-98\\
		>&{}0.
		\end{align*}
		Similarly, for the second case we have
		\begin{align*}
        HM\left(C_3\left(T^3_{n-3}\right)\right)-HM(G)
		\geq&{} HM\left(C_3\left(T^3_{n-3}\right)\right)-HM\left(C_3\left(P_3,n-5\right)\right)\\
		=&{}\left(n^3-7n^2+24n+10\right)-\left(n^3-7n^2+22n+40\right)\\
		=&{}2n-30\\
		>&{}0.
		\end{align*}
		This means that in both cases $HM(G)<HM\left(C_3\left(T^3_{n-2}\right)\right)$.
 	\item[\textit{(iii)}] $k=3$, then $G=C_3\left(T_1,T_2,T_3\right)$. By Corollary \ref{cor: cor00}, it is simple to see that $HM(G)\le HM\left(C_3\left(l_1,l_2,l_3\right)\right)$. On the other hand, since by assumption $G\not \cong C_3\left(1,1,n-5\right)$, the Hyper-Zagreb index attains its maximum when\\
 	 $G\cong C_3\left(1,2,n-6\right)$. Hence,
	\begin{align*}
	HM\left(C_3\left(T^3_{n-3}\right)\right)-HM(G)
	\geq&{} HM\left(C_3\left(T^3_{n-3}\right)\right)-HM\left(C_3\left(1,2,n-6\right)\right)\\
	=&{}\left(n^3-7n^2+24n+10\right)-\left(n^3-10n^2+43n+62\right)\\
    >&{}0.
	\end{align*}
	\end{enumerate}
\par
\textit{Case 2}: $m=4$.
	This needs to be analyzed for $k=1,2,3$ and $4$.
	\begin{enumerate}
	\item[\textit{(i)}] $k=1$, then $G=C_4\left(T_1\right)$. Since $G\not \cong C_4\left(n-4\right)$, we have $T_1\not \cong S_{n-3}$. Note that $G$ has a maximum value of the Hyper-Zagreb index if $T_1\cong T^1_{n-3}$ by Theorem \ref{thm: thm1} and Lemma \ref{lem: lem3}. Moreover, we have (for $n\ge 15$)
	\begin{align*}
	HM\left(C_3\left(T^3_{n-3}\right)\right)-HM(G)
	\geq&{} HM\left(C_3\left(T^3_{n-3}\right)\right)-HM\left(C_4\left(T^1_{n-3}\right)\right)\\
	=&{} \left(n^3-7n^2+24n+10\right)-\left(n^3-7n^2+22n+20\right)\\ 
	>&{}0.
	\end{align*}
	\item[\textit{(ii)}] $k=2$, then $G=C_4^{u_1,u_2}(T_1,T_2)_{\alpha}=C_4\left(T_1,T_2\right)_{\alpha}$, where $\alpha=d_G(u_1, u_2)$. By Lemma \ref{lem: lem4}, $G$ attains maximum $HM$-value if $G\cong C_4\left(l_1,l_2\right)_{\alpha=1}$. This lemma also implies that 
	\begin{align*}
	HM\left(C_4\left(l_1,l_2\right)\right)_{\alpha=1}\le HM\left(C_4\left(1,n-5\right)\right)_{\alpha=1}
	=n^3-7n^2+22n+38.
	\end{align*}
	Therefore, for $n\ge 15$, we have
	\begin{align*} 
	HM(G)\le n^3-7n^2+22n+38<n^3-7n^2+24n+10=HM\left(C_3\left(T^3_{n-2}\right)\right).
	\end{align*}
	\item[\textit{(iii)}] $k=3$, then $G$ is considered as $C_4\left(T_1,T_2,T_3\right)$. By Corollary \ref{cor: cor00} and Lemmas \ref{lem: lem4}, \ref{lem: lem3}, for $n \ge 15$ we have 
	\begin{align*}
	HM\left(C_3\left(T^3_{n-2}\right)\right)-HM(G)
	\geq&{} HM\left(C_3\left(T^3_{n-2}\right)\right)-HM\left(C_4\left(l_1,l_2,l_3\right)\right)\\
	>&{}HM\left(C_3\left(T^3_{n-2}\right)\right)-HM\left(C_4\left(1,n-5\right)\right)_{\alpha=1}\\
	=&{}\left(n^3-7n^2+24n+10\right)-\left(n^3-7n^2+22n+38\right)\\
	>&{}0.
	\end{align*}
	\item[\textit{(iv)}] $k=4$, then $G=C_4(T_1,T_2,T_3,T_4)$. In a similar way, one can see easily that $HM(G)<HM\left(C_3\left(T^3_{n-2}\right)\right)$; completing the proof of the second case. 
	\end{enumerate}
\par
\textit{Case 3}: $m\geq 5$.
	Using Lemmas \ref{lem: newlem}, \ref{lem: lem3} and Corollaries \ref{cor: cor0}, \ref{cor: cor1}, we conclude that (for $n\ge 15)$
	\begin{align*}
    HM\left(C_3\left(T^3_{n-2}\right)\right)-HM(G)
    \geq&{} HM\left(C_3\left(T^3_{n-2}\right)\right)-HM\left(C_m\left(l_1, l_2, \dots, l_k\right)\right)\\
    \geq&{} HM\left(C_3\left(T^3_{n-2}\right)\right)-HM\left(C_m\left(n-m\right)\right)\\
    \geq&{} HM\left(C_3\left(T^3_{n-2}\right)\right)-HM\left(C_5\left(n-5\right)\right)\\
    =&{}\left(n^3-7n^2+24n+10\right)-\left(n^3-7n^2+20n+30\right)\\
	>&{}0.
	\end{align*}
\end{proof}
\section{Conclusion}
In this paper, we  studied the  Hyper-Zagreb index and characterized the trees and unicyclic graphs with the first four and first eight greatest $HM$-value. It would be of interest to investigate its behavior on other classes of graphs with simple connectivity patterns and cyclic structures.

\providecommand{\bysame}{\leavevmode\hbox
	to3em{\hrulefill}\thinspace}


\end{document}